\input amstex

\documentstyle{amsppt}

\def\Pic{{\text{\rm Pic}}}
\def\dim{{\text{\rm dim}}}

\def\PX{{\text {\rm Pic}}^0(X)}

\def\DI{{\text {\rm dim}}   }
\def\DX{{\text {\rm dim}}(X)}

\def\alb{\text{\rm alb}}
\def\Alb{\text{\rm Alb}}
\def\lra{\longrightarrow}
\def\ot{{\otimes}}
\def\OO{{\Cal{O}}}
\def\PP{{\Cal{P}}}
\def\PPP{{\Bbb{P}}}
\def\ox{{\omega _X}}
\def\FF{{\Cal{F}}}
\def\HH{{\Cal{H}}}
\def\GG{{\Cal{G}}}
\def\BB{\Cal{B}}
\def\LL{\Cal{P}}
\def\KK{{\Cal{K}}}
\def\QQQ{{\Cal{R}}}
\def\QQ{{\Bbb{Q}}}
\def\CC{{\Bbb{C}}}
\def\sh{{\hat {\Cal     S}}}
\def\ss{{ {\Cal     S}}}

\topmatter

\title {On algebraic fiber spaces over varieties of
maximal Albanese
dimension}\endtitle
\rightheadtext{Algebraic fiber spaces }
\author Jungkai A. Chen, Christopher D. Hacon   \endauthor

\address Jungkai Alfred Chen, Department of Mathematics, National Chung
Cheng University, Ming 
Hsiung, Chia Yi, 621, Taiwan \endaddress

\email  jkchen\@math.ccu.edu.tw  \endemail

\address Christopher Derek Hacon, University of California Riverside,
Department of Mathematics, Riverside, CA 92521, USA 
\endaddress

\email   chhacon\@math.utah.edu   \endemail

\footnote""{The first author was partially supported by
National Science Council, Taiwan (NSC-89-2115-M-194-029
)}

\abstract
We study algebraic fiber spaces $f:X\lra Y$ where $Y$ is of maximal Albanese
dimension. In particular we give an effective version a theorem of
Kawamata:  If $P_m(X)=1$ for some $m \ge 2$,
then the Albanese map of $X$ is surjective. Combining this with \cite {CH}
it follows that $X$ is birational to an abelian variety
if and only if $P_2(X)=1$ and $q(X)=
\dim (X)$.
\endabstract


\endtopmatter

\document

\heading Introduction \endheading
In this paper we combine the generic vanishing theorems of \cite{GL1},
\cite {GL2}, the techniques of \cite{EL1} and the results of \cite{Ko1}
and \cite{Ko2} to answer a number of natural questions concerning
the geometry and birational invariants of irregular complex
algebraic varieties.

Throughout the paper, we are motivated by the following:

\proclaim{Conjecture K (Ueno)} Let $X$ be a
nonsingular projective algebraic variety
such that $\kappa (X)=0$ and let $\alb _X :X\lra \Alb (X)$ be the Albanese
map.
Then

(1) $\alb _X$
is surjective and has connected fibers, i.e. $\alb _X$ is an algebraic fiber
space.

(2) Let $F$ be a general fiber of $\alb _X$. Then $\kappa (F)=0$.

(3) There is an \'etale covering $B\lra \Alb (X)$ such that $X \times
_{\Alb (X)}B$ is birationally equivalent to $F\times B$ over $B$.
\endproclaim

The main evidence towards this conjecture is given by
\proclaim {Theorem (Kawamata)}

(1) Conjecture K (1) is true \cite {Ka1, Theorem 1}.

(2) If $F$ has a good minimal model
(i.e. a model with only canonical singularities,
and such that $mK\equiv 0$ for some positive
integer $m$).
Then Conjecture K is true \cite {Ka2}.
\endproclaim

We remark here that in the proofs of our statements,
we will make use \cite {Ka1, Theorem 1} in an essential way.
As a consequence of \cite {Ka1, Theorem 1}, one sees
that:
\proclaim {Corollary (Kawamata, \cite {Ka1})}
If $\kappa (X)=0$ then $q(X)\leq \dim(X)$. Moreover,
if $q(X)=\dim (X)$ then $\alb _X:X\lra \Alb(X)$ is
a birational morphism.
\endproclaim
There has also been considerable interest in
effective versions
of this result. Koll\'ar has shown that:

\proclaim {Theorem (Koll\'ar, \cite {Ko4})}

(1) If $P_3(X)=1$ then
the Albanese map is surjective.

(2) Moreover, if $P_4(X)=1$ and $q(X)=\dim(X)$,
then  $\alb _X:X\lra \Alb(X)$ is
a birational morphism.
\endproclaim

\proclaim {Conjecture (Koll\'ar, \cite {Ko3, 18.13})}
$X$ is birational to an abelian variety if and only if $q(X)=\dim (X)$
and $P_m(X)=1$ for some $m\geq 2$.
\endproclaim

Our first result is:

\proclaim {Theorem 1} If for some integer $m\geq 2$
the $m$-th plurigenera $P_m(X)$ equals $1$, then the Albanese map of
$X$ is surjective.
\endproclaim

Combining Theorem 1 with \cite {CH} we obtain:

\proclaim {Corollary 2} Koll\'ar's conjecture above \cite {Ko3, 18.13} is
true.
\endproclaim

We are also able to generalize the corollary of Kawamata's Theorem above:

\proclaim {Corollary 3}
If $P_m(X)=P_{2m}(X)=1$ for some $m\geq 2$, then
$q(X)\leq \dim (X)-\kappa (X)$. If $q(X)= \dim (X)-\kappa (X)$,
then the general fiber of the Iitaka fibration of $X$ is birationally
equivalent to a fixed \'etale cover $\tilde {A}$ of $A:=\Alb (X)$.
\endproclaim
Next, we study algebraic fiber spaces $f:X\lra Y$ where $X,Y$ are smooth
projective varieties and $Y$ is of maximal Albanese dimension.
The generic vanishing theorems of Green and Lazarsfeld are a
very effective technique in the study of irregular varieties.
In \S 2, we prove a more general version of Theorem 2.1.3 below,
which applies
to irregular varieties not necessarily of maximal Albanese dimension.
Using this result
we are able to show:

\proclaim{Theorem 4}
If $\kappa (X)=0$, let $a:=\alb _X : X \to \Alb (X)$ be the Albanese map,
then

(1) Either ${a}_* \ox$ is a zero sheaf or a torsion line bundle.

(2) $P_1(F_{X/\Alb (X)})\leq 1$.

(3) If $P_1(X)=1$, then $a_* \ox = \OO_{\Alb (X)}$.

(4) There is a generically finite cover $\tilde {X}\lra X$ with $\kappa
(X)=\kappa (\tilde {X})$, such that
${\alb _{\tilde {X}}}_* (\omega _{\tilde {X}})
=\OO _{\Alb_{\tilde {X}}}.$
\endproclaim

We remark that the statement $(2)$
can be regarded as an effective version of Conjecture K (2).
Using a well known result of Fujita \cite {Mo, (4.1)},
for any $X$ with $\kappa (X)=0$,
there exists a generically finite cover $\tilde {X}\lra X$ such that
$\kappa (\tilde {X})=0$ and $P_1(\tilde {X})=1$ as in (3).

We explain the significance of (4).
For any algebraic fiber space $f:X\lra Y$
the rank at a generic point of the sheaves $f_*(\omega _{X/Y}^{\ot m})$
corresponds to the plurigenera $P_m(F_{X/Y})$ of a generic geometric fiber.
It is expected that the positivity of the sheaves $f_*(\omega _{X/Y}^{\ot
m})$
measures the birational variation
of the geometric fibers.
Our methods unluckily cannot be applied to the case $m\geq 2$
because of the lack of a suitable geometric interpretation of
the sheaves ${f}_* (\omega ^{\ot m}_{{X}/Y})$.

If $F=F_{X/\Alb (X)}$ has a good minimal model, then by Kawamata's
result, there is an \'etale covering $B\lra \Alb (X)$ such that
$X\times _{\Alb (X)}B$ is birationally equivalent to $F\times B$ over $B$.
Since $\kappa (F)=0$, by Fujita's lemma again, there exists
a generically finite cover $\tilde {F}$ such that $\kappa (\tilde {F})=0$
and $P_1(\tilde {F})=1$.
Since $F$ has a good minimal model, we may assume in fact
that (for an appropriate birational model with canonical singularities)
$K_{\tilde {F}}=0$.
Let $\tilde {X}=\tilde {F}\times B$. Then, for all $m\geq 1$
$${\pi _{B}}_* (\omega ^{\ot m}_{\tilde {X}})
=\OO _{B},$$ where $\pi _B:\tilde {X}\lra B$ is the projection
to the second factor.

The following lemma is useful:

\proclaim{Lemma 5} Let $f:X\lra Y$ be an algebraic fiber space,
$Y$ of maximal Albanese dimension. If $\kappa (X)=0$, then $Y$ is birational
to an abelian variety.
\endproclaim

\proclaim{Theorem 6}
Let $f:X\lra Y$ be an algebraic fiber space, $Y$ of maximal Albanese
dimension.

(1) If $\kappa (X)=\kappa (Y)$, then $P_1(F_{X/Y})\leq 1$.


(2) If $P_1(F_{X/Y})\ge 1$, then $ \kappa (X)\geq \kappa (Y)$.
\endproclaim

The above statements are closely connected to the following
well known conjecture

\proclaim {Conjecture $C_{n,m}$}
Let $f:X\lra Y$ be an algebraic fiber space,
$\dim (X)=n$ and $\dim(Y)=m$. Then
$$\kappa (X)\geq \kappa (Y)+ \kappa (F_{X/Y}).$$
\endproclaim

This conjecture is true when $F_{X/Y}$ has a good minimal model \cite {Ka2}.
If one could generalize the generic vanishing theorem to the sheaves
of the form $\omega _Y \ot f_* \omega _{X/Y} ^m$ for $m\geq 2$, then
using the same techniques,
the $C_{n,m}$ conjecture would follow for all algebraic fiber spaces
$ f:X\lra Y$, with $Y$ of maximal Albanese dimension.

Finally, we prove a generalization of \cite {Ka1, Theorem 15}:

\proclaim{Theorem 7}
Suppose $X$ is a variety with  $\kappa (X)=0$ and $\dim (X)\leq 2q(X)$.
If $P_1(X)>0$, then $h^0(X,\Omega _X^{n-q})>h^0(\Alb,\Omega _{\Alb}
^{n-q})$.
In particular, if $\dim (X)=q(X)+1$, then $P_1(X)=0$.
\endproclaim
We then illustrate how one can use this result to recover
the Conjecture K, in the case $\dim (X)=q(X)+1$ \cite {Ka1, Theorem 15}.
\bigskip
\noindent {\bf Acknowledgment.}
We are in debt to P. Belkale, A. Bertram, L. Ein,
J. Koll\'ar and R. Lazarsfeld for valuable coversations.

\heading 0. Conventions and Notations   \endheading

(0.1) Throughout this paper, we work
over the field of complex numbers $\CC$. $X$ will denote a smooth
projective algebraic variety.

(0.2) A general point of $X$ denotes a point in the complement
of a countable union of proper Zariski closed subsets of $X$.

(0.3) Let $f:X\lra Y$ be an algebraic fiber space, i.e. a surjective
morphism
with connected fibers. Then $F_{X/Y}$ will denote the fiber over a general
geometric point of $Y$.

(0.4) For $D_1,D_2$ $\QQ$-divisors on a variety $X$,
we write $D_1 \equiv D_2$ if $D_1$ and $D_2$ are numerically equivalent.

(0.5) For a real number $a$, let $\lfloor a \rfloor$ be the
largest integer $\le a$
and $\lceil a \rceil$ be the smallest integer $\ge a$.
For a  $\QQ$-divisor $D=\sum a_i D_i$, let
$\lfloor D \rfloor =\sum \lfloor a_i \rfloor D_i$ and
$\lceil D \rceil =\sum \lceil a_i \rceil D_i$.
 We will say that a $\QQ$-divisor $\Delta$ is $klt$ (Kawamata log terminal)
if $\Delta $ has normal crossings support and $\lfloor \Delta \rfloor =0$.
We refer to \cite {Ko3, 10.1.5} for the general
definition of klt divisor.

(0.6) Let $\FF$ be a coherent sheaf on $X$,
then $h^i(X,\FF)$ denotes the complex dimension of $H^i(X,\FF)$.
In particular,
the plurigenera $h^0(X,\ox^{ \otimes m})$ are denoted  by $P_m(X)$
and the irregularity $h^0(X,\Omega_X^1)$ is denoted by $q(X)$.

\head 1. Preliminaries \endhead

\subhead 1.1 The Iitaka Fibration \endsubhead

Let $X$ be a smooth complex projective variety with $\kappa (X)>0$. Then
a nonsingular representative of the Iitaka fibering of $X$ is a morphism
of smooth complex projective varieties $f':X'\lra V$ such that
$X'$ is birational to $X$,
$\dim (V)=\kappa (X)$ and $\kappa ({X'}\hskip-.1cm _v)=0$,
where ${X'}\hskip-.1cm _v$ is a generic geometric fiber of $f'$.
Since our questions will be birational in nature,
we may always assume that $X=X'$.
Let $A:=\Alb (X)$ and let $Z$ denote the image of $X$ in $A$. Let $Z'$
denote
an appropriate desingularization of $Z$, we may assume that $X\lra Z$ factor
s
through $Z'$.
By \cite {Ka1} the images $a(X_v)=K_v$ are translates of abelian
subvarieties
of $A$. Since $A$ contains at most countably many abelian subvarieties,
we may assume that $K_v$ are all translates of a fixed abelian subvariety
$K\subset A$. Let $S:=A/K$ and $W$ denote the image of $Z$ in $S$.
Let $W'$ be an appropriate
desingularization of $W$. We may assume that the induced morphism
$\pi :X\lra W$ factors through a morphism $\pi ' : X\lra W'$.
Consider now a birational model of the Iitaka fibration $f:X\lra V$.

{\bf Claim.}
We may assume that the map $\pi '$ factors through $f$ and a morphism
$q':V\lra W'$.

To see this, note that by construction, there is an open dense subset $U$ of
$V$ and a map $U\lra S$. However by a standard argument this must complete
to a rational map $V\lra S$
(see e.g. \cite {Ka1, Lemma 14}).
Since the problem is birational, we may assume that $V\lra S$
is infact a morphism that factors
through $W'$. The above maps fit in the following comutative diagram

$$
\CD
X &@>>>&Z'&@>>>&Z&@>{\subset}>>&A \\
@V{f}VV & & & & & & & & & & @VV{p}V \\
V & @>{q'}>>&W'&@>>>&W&@>{\subset}>>& S.\\
\endCD
$$

\bigskip

\subhead 1.2 Fourier Mukai Transforms \endsubhead

Let $A$ be an abelian variety, and denote the corresponding dual
abelian variety by $\hat {A}$.
Let $\PP$ be the normalized Poincar\'{e} bundle on
$A\times \hat {A}$. For any point $y\in \hat {A}$ let
$\PP _y$ denote the associated topological trivial line bundle.
Define the functor $\sh $ of $\OO_A$-modules
into the category of $\OO_{\hat {A}}$-modules by
$$\sh (M)=\pi _{\hat {A},*}(\PP \ot \pi ^*_AM).$$
The derived functor $R\sh$ of $\sh$ then induces an equivalence of
categories between the two derived categories $D(A) $ and $D(\hat {A})$.
In fact, by \cite{M}: {\it There are isomorphisms of functors:
$$R\ss \circ R\sh \cong (-1_A)^*[-g]$$
and
$$R\sh \circ R\ss \cong (-1_{\hat {A}})^*[-g],$$
where $[-g]$ denotes "shift the complex $g$ places to the right".}

The index theorem (I.T.) is said to hold for a coherent sheaf $\FF$ on $A$
if there exists an integer $i(\FF )$ such that for all $j
\ne i(\FF )$,
$H^j(A,\FF \ot P)=0$ for all $P\in \Pic ^0(A)$.
The weak index theorem (W.I.T.) holds for a coherent sheaf $\FF$ if
there exists an integer which we again denote by $i(\FF )$ such that for all
$j
\ne i(\FF )$, $R^j\sh
(\FF )=0$. It is easily seen that the I.T. implies the W.I.T.
We will denote the coherent
sheaf $R^{i(\FF )}\sh (\FF )$ on $\hat {A} $ by $\hat {\FF }$.

It follows, for example that given any coherent sheaf $\FF$, if
$h^i(A,\FF \ot P)=0$ for all $i$ and all $P\in \Pic ^0(A)$ then
$\FF =0$.
Consequently, if $\FF \lra \GG$ is an injection of sheaves
which induces isomorphisms in cohomology
$$H^i(A,\FF \ot P)@>{\cong }>> H^i(A,\GG \ot P),$$
for all $i$ and all $P\in \Pic ^0(A)$, then $\FF \cong \GG$.

\proclaim {Proposition 1.2.1}
If $\FF $ is a coherent sheaf on $A$ such that for all $P\in \Pic ^0(A)$
we have $h^0(A,\FF \ot P)=1$ and $h^i(A,\FF \ot P)=0$ for all $i>0$.
Then $\FF $ is supported on an abelian subvariety of $A$.
\endproclaim
\demo {Proof}
$\FF$
satisfies the WIT and $\hat {\FF }$ is a
line bundle $M$ on $\Pic ^0 (A)=\hat {A}$
such that $M$ has index $i(M)=\dim (A)$, and $\hat {M }=(-1_A)^*\FF$.
Any line bundle with $i(M)=\dim (A)$, is negative semidefinite.
It is well known that there exists a morphism of abelian varieties
$b:\Pic ^0 (A)\lra A'$ such that $M=b^*M'$ for some negative definite
line bundle $M'$ on $A'$.
It follows that $\hat {M }$ and hence $\FF$ are supported on
the image of $b^*:\hat{A'}\lra A$.
\hfill $\square$
\enddemo

\head 2. Relative Generic Vanishing Theorems \endhead
We start by recalling some facts on cohomological support loci.
Let $\pi:X\to A$ be a morphism from a smooth projective variety $X$ to
an abelian variety $A$.
If $\FF$ is a coherent sheaf on $X$,
then one can define the {\it cohomological support loci} by

$$V^i(X,A,\FF):=\{ P \in \Pic ^0(A) | h^i(X, \FF \otimes \pi^*P) \ne 0
\}.$$
In particular,
if $\pi={ {\alb }}_X :X \to \Alb (X)$,
then we simply write
$$V^i(X,\FF):=\{ P \in \Pic ^0(X) | h^i(X, \FF \otimes P) \ne 0 \}.$$

We say that $X$ has maximal Albanese dimension if $\DI {(\alb _X(X))}=\DI
(X)$.
The geometry of the loci $V^i(X,\ox )$ defined above
is governed by the following:

\proclaim{Theorem 2.1 (Generic Vanishing Theorem)}

(1) Any irreducible component of $V^i(X,\ox )$ is a translate of
a sub-torus of $\Pic ^0(X)$
and is of codimension at least
$i-(\DI (X) -\DI (\alb _X(X)))$.

(2) Let $P \in T$ be a general point of an irreducible component
$T$ of $V^i(X,\ox)$.
Suppose that $v
\in H^1(X, \OO
_X)\cong T_{P }\Pic ^0(X)$ is not tangent to $T$. Then the sequence
$$H^{i-1}(X,\ox \otimes P) \overset {\cup v} \to {\longrightarrow}
H^{i}(X,\ox \otimes P)  \overset {\cup v} \to {\longrightarrow }
H^{i+1}(X, \ox \otimes P) $$
is exact.
If $v$ is tangent to $T$, then the maps in the above sequence vanish.

(3) If $X$ is a variety of maximal Albanese dimension, then
$$\Pic ^0(X)\supset V^0(X,\ox )\supset V^1(X,\ox)
\supset \ ...\ \supset V^n(X,\ox )=\{ \OO _X\}.$$

(4) Every irreducible component of $V^i(X,\ox)$ is
a translate of a sub-torus of $\PX$ by a torsion point.
\endproclaim

\demo {Proof} For (1),(2), see \cite{GL1},\cite{GL2}. For (3), see
\cite{EL1} and
for (4), see \cite{S}.
\hfill $\square$ \enddemo

In \cite{EL1}, Ein and Lazarsfeld illustrate various examples
in which the geometry of $X$ can
be recovered from information on the loci $V^i(X, \ox)$.

In this section we prove a relative version of Theorem 2.1.
Let $\pi :X\lra Y$ be a surjective map of smooth projective varieties.
Assume that $Y$ has maximal Albanese dimension.
Let $n$ and $m$ be the dimension of $X$ and $Y$ respectively.
We wish to study the geometry of the loci $V^{i}(Y,R^j\pi _*
\omega _X )$. By a result of Simpson \cite{S}, the irreducible
components of these loci are torsion translates of subtori
of $\Pic ^0(Y)$. Therefore, their geometry is completely determined
by their torsion points. Recall \cite{Ko2, Corollary 3.3} that for
any torsion $Q\in \Pic ^0(X)$, one has $$R^{\cdot }\pi _* (\omega
_X \ot Q)\cong \sum R^if_*(\omega _X \ot Q)[-i].$$
In particular $h^p(X,\omega _X\ot Q)=\sum h^i(Y,R^{p-i}f_*(\omega _X \ot
Q))$.

\proclaim{Proposition 2.2} For $\pi :X\lra Y$ as above
$$\Pic ^0(Y)\supset
V^0(Y,R^j\pi _* \omega _X)\supset V^{1}(Y,R^j\pi _*
\omega _X )\supset $$
$$\hskip3cm ... \supset
V^m(Y,R^j\pi _* \omega _X ).$$
\endproclaim

\demo{Proof} We will prove that $h^i(Y,R^j\pi _* \omega _X\ot P)>0$ implies
$h^{i-1}(Y,R^j\pi _* \omega _X\ot P)>0$ for all $i>0$.
By Simpson's result mentioned above, it suffices to prove the assertion
for torsion elements
$P\in \Pic ^0(Y)$.
Fix a very ample line bundle $L$ on $Y$ and for $1\leq \alpha \leq m$,
let $D_{\alpha}$ be
sufficiently general divisors in $|L|$. Denote the preimage by
$\tilde {D}_{\alpha}= \pi ^{-1}(D_{\alpha})$.
Consider the sequence of varieties defined by $X_0=X$, $X_{\alpha +1}=
X_{\alpha}\cap
\tilde {D}_{\alpha +1}$. Let $Y_{\alpha}:=\pi (X_{\alpha})$.
We may assume that the $D_{\alpha}$ are chosen so that
$X_{\alpha}, Y_{\alpha}$ are smooth for $0\leq {\alpha} \leq m$.
The corresponding morphisms $\pi |_{X_{\alpha}}:X_{\alpha}
@>>>Y_{\alpha }$ will also be denoted simply by $\pi$.

By \cite{Ko 2, Corollary 3.3}, for any torsion $P\in \Pic ^0(Y)$, we may
identify
$H^i(Y_{\alpha},R^j$ $\pi _* \omega _{X_{\alpha}}\ot P)$
as a subgroup of $ H^{i+j}(X_{\alpha},$ $\omega _{X_{\alpha}}\ot  \pi ^*P)$.
\proclaim{Claim}
There exists a surjective map
$$H^0(Y_i,R^j\pi _*\omega _{X_i}\ot P)\lra H^i(Y,R^j\pi_* \omega _X\ot P).$$
\endproclaim

\demo{Proof} We will prove the assertion for $P=\OO _Y$. In the general
case the proof proceeds analogously. Consider the exact sequence of sheaves
$$0\lra \omega _{X_t}\lra \omega _{X_t}(\tilde{D}_{t+1})\lra
\omega _{X_{t+1}}\lra 0$$
this induces sequences of sheaves
$$0\lra R^j\pi _*\omega _{X_t}\lra R^j\pi _*
\omega _{X_t}(\tilde {D}_{t+1})\lra
R^j\pi _*\omega _{X_{t+1}}\lra 0.$$

By step 4 of the proof of \cite {Ko1, Theorem 2.1 iii},
for appropriately chosen $D_{\alpha}$,
the above sequence is exact and equivalent to
$$R^j\pi _*\omega _{X_t} \ot \left( 0\lra \OO _{ {Y}_t}\lra \OO  _{{Y}
_t}({D}_{t+1})
\lra \OO  _{{Y}_{t+1}}({D}_{t+1})\lra 0\right) .$$

By \cite{Ko 1}, $H^k({Y}_t,R^j\pi _*\omega _{X_t}({D}_{t+1}))=0$
for all $k>0$.
We therefore have an exact sequence
$$0\lra H^0({Y}_t,R^j\pi _*\omega _{X_{t}})\lra H^0({Y}_t,R^j\pi _*
\omega _{X_{t}}({D}_{t+1}))\lra $$
$$\ \ \ \ \ \ \ \ \ \ \ \ \ \ \ \ \ \ \ \ \ \ \ \
H^0({Y}_{t+1},R^j\pi _*\omega _{{X}_{t+1}})\lra
H^1(Y_t, R^j\pi _*\omega _{X_{t}})\lra 0,$$
and isomorphisms
$H^k({Y}_{t+1},R^j  \pi _* \omega _{X_{t+1}}) \cong H^{k+1}
({Y}_t,R^j\pi _* \omega _{X_t})$ for $k=1,...,m-t-1.$
The claim now follows.
\enddemo

Moreover, the map $$H^0({Y}_{i},R^j\pi _* \omega _{X_i})\lra H^i(Y,
R^j\pi _* \omega _{X})\lra H^{i+j}(X,\omega _X)$$
is induced by the inclusion $H^0(Y_{i},
R^j \pi _* \omega _{X_i})\lra H^j(X_i,\omega
_{X_i})$ and by the coboundary maps
$$\delta :H^j(X_i,\omega _{X_i})@>{\delta _{i-1}^j}>> H^{j+1}(X_{i-1},
\omega _{X_{i-1}})
@>{\delta _{i-2}^{j+1}}>> ... @>{\delta _{0}^{j+i}}>>
H^{i+j}(X,\omega _X ).$$

The map $\delta$ is  dual to the map
$$\delta ^*:H^{n-i-j}(X,\OO _X)\lra H^{n-i-j}(X_1,\OO _{X_1})\lra ...
\lra H^{n-i-j}(X_i,\OO _{X_i})$$
induced by successive restrictions.
In turn, this map is complex conjugate to the map
$$\bar {\delta ^*}:H^{0}(X,\Omega  ^{n-i-j}_X)\lra H^{0}
(X_1,\Omega ^{n-i-j}_{X_1})\lra ...
\lra H^{0}(X_i,\Omega^{n-i-j} _{X_i})$$
induced by successive restrictions.
Let $V$ be the subspace of $H^{0}(X,\Omega  ^{n-i-j}_X)$
corresponding to the complex conjugate of ${H^i(Y,R^j \pi _* \omega _X)^*}$.
It follows from the above claim that $\bar {\delta }^*$ induces
an injection $V\hookrightarrow H^{0}(X_i,\Omega^{n-i-j} _{X_i})$.

Fix a general smooth point $p$ of $\pi (X_i)$.
Let $a: Y \lra A$ be the Albanese map.
For a general point $p \in Y_i$.
One may assume that $z_1,...,z_m$ are local
holomorphic coordinates for $a(Y)$ at $a(p)$, where $z_1,...,z_{m-i}$ are
local
holomorphic coordinates for $a(Y_i)$ at $a(p)$.
We may assume furthermore that the pull backs of the $z_i$, which
we denote by $x_i$, give rise to  local
holomorphic coordinates for $Y$ and $Y_i$ at $p$.
We may assume that $dx_l(p)=0$ in $\Omega ^1 _{Y_i}\ot \CC (p)$
for all $m-i+1\leq l\leq m$ (i.e. $dx_l$ is conormal to $Y_i$ at $p$).
Fix $\tilde {p}\in \pi ^{-1} (p)$.
Let $\omega \in H^0(A,\Omega ^1 _A)$
such that $\omega (p)=
dz_{m}$.
Let $y_1,...,y_n$ be local
holomorphic coordinates on $X$ centered at $\tilde {p}\in \pi ^{-1}(p)$
such that $y_i=x_i \circ \pi$ for all
$1\leq i\leq m$. Let $v$ be any element in
$V\subset H^0(X,\Omega ^{n-j-i}_X)$.
Then,
$$v(\tilde {p})=\sum a_Kdy_K$$
where the sum runs over all multiindices of length $n-i-j$ ie
$K=(k_1,...,k_{n-i-j})$ and $dy_K=dy_{k_1}\wedge ...\wedge dy_{k_{n-i-j}}$.
Since $V$ injects in $H^0(X_i,\Omega _{X_i}^{n-i-j})$,
for any fixed $v\in V$, by genericity of the choice of the point $p$,
we may assume that $v(\tilde {p}) \ne 0$ in
$\Omega _{X_i}^{n-i-j} \ot \CC
(\tilde {p})$. In particular,
there exists a multiindex $\bar {K}$ such that $a_{\bar {K}}(\tilde {p})
\ne 0$, and for all $1\leq l\leq n-i-j$, $\bar {k}_l $
doesn't belong to $\{ m-i+1,...,m\}$. Therefore $v\cup \omega \in
H^0(X,\Omega ^{n-j+1}_X)$ is non-zero since the
coefficient of $dy_{ \bar {K}\cup \{ m\} }$ is non-zero.

Composing again with complex conjugation and Serre duality
we see that there
is a non-zero element in $H^{i-1} (Y,R^j \pi _* \omega _X)
\subset H^{j+i-1}(X,\omega _X)$.
\hfill $\square$ \enddemo

The following notation will be convenient.
For any line bundle $L$ on $Y$ and $v\in H^1(Y,\OO _Y)$,
we will denote by $\KK ^i_{L,Y,v}$ and $\BB ^{i+1}_{L,Y,v}$ the kernel and
the image
of the map
$$H^i(Y,L)@>{\wedge v}>>H^{i+1}(Y,L).$$
Let $\HH ^{i}_{L,Y,v}:=\KK ^i_{L,Y,v}/\BB ^{i}_{L,Y,v}$.
The subscripts $Y$ and $v$ will be dropped when no confusion is likely.

Let $\tau :H^1(X,\OO _X )@>>> \Pic ^0(X)$ be the map induced by the
exponential sheaf sequence. Let $ \Delta \subset Spec \CC [t]$
be a neigborhood
of $0$. For $P\in \Pic ^0(X)$ and $v\in H^1(X,\OO _X )$, let
$\LL$ a line bundle on $X\times \Delta$ such that $\LL |_{X\times t}
\cong P\ot \tau (tv)$. Let $p_X:X\times \Delta @>>>X$ and
$p_{\Delta}:X\times \Delta @>>>\Delta$ be
the projections to the first and second factors. We will need the following
result due to Green and Lazarsfeld.

\proclaim {Theorem 2.3 \cite{GL2}} There is a neighborhood of $0$ for which
$$R ^i{p_{\Delta_*}}
(p_X^*\omega _X \ot \LL ) \cong \left( \KK _{\omega _X\ot P}^i
\ot \OO _{\Delta / m}\right) \oplus
\left( \HH _{\omega _X \ot P}^i \ot \OO
_{\Delta}\right) .$$
\endproclaim
\demo{Proof} This is a generalization of \cite{GL2, Theorem 3.2},
which follows from the comments preceeding \cite{GL2, Theorem 6.1}.
See also \cite {ClH}.
\hfill $\square$ \enddemo

\proclaim {Corollary 2.4} Let $\phi
\in \KK _{\omega _X \ot P}^i \subset H^i(X,\omega _X \ot
P)$. Let $\Phi \in R ^i{p_{\Delta_*}}
(p_X^*\omega _X\ot \LL ) $ be a
section  such that
$(\Phi)_0=\phi$. Assume that
$\Phi |_{(\Delta -0 )}=0$.  Then $\phi =\gamma \cup v$ for an appropriate
$\gamma \in H^{i-1}(X, \omega _X \ot P)$.
\endproclaim

\demo{Proof} Since $\Phi |_{(\Delta -0 )}=0$, there exists an integer
$k\geq 0$ such that $t^k\Phi =0\in R ^i{p_{\Delta_*}} (p_X^*\omega _X \ot
\LL ) $.
Therefore, $$\phi t^k +\phi _1t^{k+1}+\phi _2 t^{k+2}+\ldots =0 \in
\left( \KK _{\omega _X \ot P} ^i \ot \OO _{\Delta / m}\right) \oplus
\left( \HH _{\omega _X \ot P}^i \ot \OO _{\Delta}\right) .$$
If $k=0$ then $\phi =0 \in \KK _{\omega _X \ot
P}^i$, and there is nothing to show.
If $k\geq 1$ then $\phi =0 \in \HH _{\omega _X \ot P}^i$, and hence $\phi  =
\gamma \cup v$ for an appropriate $\gamma \in H^{i-1}(X, \omega _X \ot P)$.
\hfill $\square$ \enddemo

\proclaim {Proposition 2.5} Let $\pi:X@>>> Y$ be an algebraic
fiber space, $Y$ of maximal Albanese dimension.
If $H^i(Y,\pi _*\omega _X\ot P)=0$ for all $i$ and
all $P$ in a punctured neighborhood of a torsion $P_0 \in \Pic ^0(Y)$,
then, for any $v\in H^1(Y,\OO _Y)$, the complex
$$H^{i-1} (Y,\pi _*\omega _X \ot P_0)
@>{\cup v}>>H^{i} (Y,\pi _*\omega _X \ot P_0)@>{\cup v}>>
H^{i+1} (Y,\pi _*\omega _X \ot P_0)\ \ \ \ \ \ \ \ \ \ \ \ (*)$$
is exact.
\endproclaim

\demo{Proof}
We will use the notation of the proof of Proposition 2.2.
Let us first consider the folowing diagram, which is commutative for each
square.
$$
\CD
{p_{\Delta }}_* (p_X^*\omega _{X_i} \ot \LL )  @>{\delta_{\Delta}}>>
R^i {p_{\Delta_* }}(p_X^*\omega _X \ot \LL ) \\
@V{res_{t=0}}VV  @V{res_{t=0}}VV\\
H^0(X_i,\omega _{X_i}\ot \pi ^* P_0) @>{\delta_{X}}>> H^i(X,\omega _X \ot
\pi ^* P_0)
  @>{\cup \pi^*v}>>  H^{i+1}(X,\omega _X \ot \pi ^* P_0) \\
@A{\parallel}AA  @A{\cup}AA   @A{\cup}AA \\
H^0(Y_i,\pi _* \omega _{X_i}\ot P_0)  @>{\delta_{Y}}>>  H^i(Y,\pi _* \omega
_{X_i}\ot P_0)
 @>{\cup v}>>  H^{i+1}(Y,\pi _* \omega _{X_i}\ot P_0)
\endCD
$$

Let $f\in H^i(Y,\pi _* \omega _X \ot P_0)$, such that $f\cup v =0 \in
H^{i+1}(Y,\pi _* \omega _X \ot P_0)$.
By abuse of notation, we will also denote by $f$
the corresponding element in $H^i(X,\omega _X\otimes \pi ^* P_0)$.
Since $\delta_X$ is surjective, we take $\tilde{f} \in H^0(X_i,\omega
_{X_i}\ot \pi ^* P_0)$
a lift of $f$.
We remark that also $\tilde {f}
\cup \pi ^* v=0$.
 Let  $\tilde {F} \in {{p_{\Delta }}}_*
(p_X^*\omega _{X_i} \ot \LL )$ be the section
corresponding to $\tilde{f} +0t+0t^2+...$ and
$F \in R^i{p_{\Delta_* }}(p_X^*\omega _X \ot \LL )$ be
the corresponding section under the coboundary map.
One has $F |_{t=0}=f \in H^i(X,\omega _X \ot \pi ^*P_0)$

By the claim in the proof of 2.2,
the condition $H^i(Y,\pi _*\omega _X\ot P)=0$ is equivalent to
the vanishing of the following map
$$H^0(Y_i,\pi _* \omega _{X_i}\ot P)\cong
H^0(X_i,\omega _{X_i}\ot \pi ^*P)@>>> H^i(X,
\omega _X \ot \pi ^*P)\hskip1cm (**).$$

For all $P$ in a punctured neighborhood of $P_0$ the map $(**)$ vanishes.
Therefore $F |_{(\Delta -0)}=0$. By Corollary 2.4, $f =\gamma \cup \pi^*v$
for an appropriate $\gamma \in H^{i-1}(X, \omega _X \ot \pi ^*$ $ P_0)$.

Following \cite{Ko2}, write $\gamma =\sum g_j$ where $g_j\in
H^{i-1-j}(Y,R^j\pi _*\omega _X\ot P_0)$,
then $g_j\cup v\in H^{i-j}(Y,R^j \pi _*$ $\omega _X \ot P_0^*)$,
so $f=g_0\cup v$, and therefore $(*)$ is exact.
\hfill $\square$
\enddemo

\head 3. Proof of Theorems \endhead

\demo{Proof of Theorem 1}
By \cite {Ka1, Theorem 1}, we may assume that $\kappa (X)>0$.
Let $f:X\lra V$ be a birational model of the Iitaka fibration,
$\pi  :X\lra W$ and $S$ be as in \S 1.1.

Assume that $m=2$ (the proof proceeds analogously for any $m\geq 2$).
Fix $H$ an ample divisor on $S$.
For a fixed $r \gg 0$,
after replacing $X$ by an
appropriate birational model,
we may assume that
$$|rK_X-\pi ^* H|=|M_r|+F_r$$
where $|M_r|$ is non-empty and free,
and $F_r$ has simple normal crossings.
Let $B$ be a general divisor of the linear series
$|rK_X-\pi ^* H|$. We may assume again that $B$ has normal crossing support.
Define $L:=\OO _X(K_X-\lfloor  \frac{B}{r} \rfloor )$.
We have that $L\equiv (\pi ^* H/r)+\{ \frac{B}{r} \} $ is
numerically equivalent to the sum of the pull back of
a nef and big $\QQ $-divisor on $W$ and
a klt $\QQ $-divisor on $X$.

As in the proof of  \cite {CH, Lemma 2.1}, it
is possibile to arrange that $|2K_X|=|K_X+L|+\lfloor \frac{B}{r} \rfloor $.
In particular $h^0(X,\omega _X\ot L)=1$.

Since $\pi :X\lra W$ is a surjective map, by \cite{Ko3, Corollary 10.15},
$\pi  _* (\ox \ot L)$ is a torsion free coherent sheaf on $W$,
and $h^i(W, \pi  _* (\ox \ot L)\ot P)=0$ for all $i>0$ and
$P\in \Pic ^0(S)$. It follows that
$$h^0(W, \pi _* (\ox \ot L)\ot P)=\chi (W,\pi _* (\ox \ot L)\ot P)$$
$$=\chi (W,\pi _* (\ox \ot L))=h^0(W, \pi _* (\ox \ot L))=1$$
for all $P\in \Pic ^0(S)$.

By Proposition 1.2.1, $\pi  _* (\ox \ot L)$ is supported
on an abelian subvariety $S'$ of $S$. Since $\pi  _* (\ox \ot L)$
is torsion free and the image of $X$ generates
$S$, we see that $S'=S$, and hence $X\lra A$ is surjective. \hfill $\square$
\enddemo

\demo {Proof of Corollary 3}
If $\kappa (X)=0$, this is a result of Kawamata.
We may therefore assume that $\kappa (X)>0$.
By Theorem 1, $X\lra \Alb (X) $ is surjective.
By \cite {CH, Lemma 2.1} $h^0(X, \omega _X ^{\ot m}\ot $ $\pi ^*P)>0$
for all $P\in \Pic ^0(S)$.
{From} the map of linear series
$$|mK_X + Q|\times |mK_X - Q|\lra |2mK_X|$$
it follows that $P_{2m}(X)\geq \dim (S)+1$, and therefore
$\dim (S)=0$. The general geometric fiber $X_v$
of the Iitaka fibration has
dimension $\dim (X)-\kappa (X)$. Since $\dim (S)=0$, by construction, it
follows that $\alb _X (X_v)=\Alb (X)$ and hence $\dim (X)-
\kappa (X)\geq q(X)$.

Suppose now that $\dim (X)-\kappa (X)=q(X)$.
Then the map $X_v\lra \Alb (X)$ is birationally \'etale,
and we may assume that $X_v$ is birational to $\tilde {A}$ for a fixed
abelian variety
$\tilde {A}$.
\hfill $\square$
\enddemo

\noindent {\bf Remark.}
Under the same hypothesis one can show that in fact $q(V)=0$.
In particular if $P_m(X)=P_{2m}(X)=1$ for some $m\geq 2$ and $q(X)=\dim (X)
-1$, $\kappa (X)=1$, then $V=\PPP ^1$.

\proclaim {Lemma 3.1} If $\kappa (X)=0$, then $V^0(X, \omega _X)$
contains at most one point.
\endproclaim
\demo{Proof} Assume that there are $2$ distinct points $P,\ Q$
in $V^0(X,\omega _X)$. By Theorem 2.1, we may assume that
$P,Q\in \PX$ are torsion elements. Pick any $m>0$ such that
$P^{\ot m}=Q^{\ot m}=\OO _X$, then
if $P\ne Q$, we have
$h^0(X, \omega_X^{\otimes m})>1$ which is impossible. Therefore $P = Q$.
\hfill $\square$ \enddemo

\proclaim {Lemma 3.2} Let $a:X \lra A$ be an algebraic fiber space from
a smooth projective variety $X$ to an abelian variety $A$.
Then $$K:=ker \left( \Pic ^0(A)@>>> \Pic ^0(X)\right) =\{ \OO _A \}.$$
\endproclaim
\demo{Proof}
The map $a$ factors as $X@>{\alb _X}>> \Alb (X)@>>> A$.
The map $\Alb (X)\lra A$ is surjective,
therefore $K$ consists of a discrete
subgroup of $\Pic ^0(A)$. Define $A'=(\Pic ^0(A)/K)^*$. Then there are maps
$$X\lra \Alb (X) \lra A'\lra A.$$
Since $X\lra A$ has connected fibers, then $A=A'$.
\hfill $\square$ \enddemo

\demo {Proof of Theorem 4}
(1) Let $A:=\Alb (X)$.
By \cite{Mu} $a _* \ox $ is zero if and only if $V^i(A,a _* \ox )
$ is empty for all $i$.
By Proposition 2.2, this is equivalent to $V^0(A,a _* \ox )$ being empty.
Thus if $a _* \ox  \ne 0$, by Lemma 3.1, we may assume that
$V^0(A, a _* \ox )$ consists of exactly one (torsion) point say $P$ and
$h^0(A, a _* \ox \ot P)=1$.
We shall prove that the injection $\OO _A \lra a _* \ox \ot P$
is in fact an isomorphism of sheaves. Therefore $a _* \ox \in \Pic ^0(A)$.

To this end,
we consider  complex $D(v)$
$$ .\ .\ .\ H^{i-1}(A,a_* \omega _X \ot P)
\overset {\cup v} \to {\longrightarrow}
H^{i}(A,a_* \omega _X \ot P)  \overset {\cup v} \to {\longrightarrow }
H^{i+1}(A,a_* \omega _X \ot P) \ .\ .\ . $$
By Proposition 2.5, this is exact for all $v\in H^1(A,\OO _A)$.

\noindent {\bf Step 1.} Let $V^*=H^1(A,\OO _A)$ and $\PPP =\PPP (V)$.
There is an exact sequence of vector bundles on $\PPP$:
$$0\lra H^0 (A,a_* \omega _X \ot P)\ot \OO _{\PPP} (-q)\lra
H^1 (A,a_* \omega _X \ot P)\ot \OO _{\PPP} (-q+1)\lra ...$$
$$\ \ \ \ \ \ \ \ \ \ \ \ \ \ \ \ \ \ ...
\lra H^q (A,a_* \omega _X \ot P)\ot \OO _{\PPP} \lra 0.\ \ \ \ \ \ \ \ \
(\KK ^{\bullet})$$

To see this, it is enough to check exactness on each fiber (see \cite {EL1}
for a similar argument).
A point in $\PPP$ corresponds to a line in $H^1(A,\OO _A)$
containing a point say $v$. On the fibers above $[v]$, the sequence
of vector bundles corresponds to the complex $D(v)$ which is exact.
Similarly, there is an exact sequence of vector bundles on $\PPP$:
$$0\lra H^0 (A,\OO _A)\ot \OO _{\PPP} (-q)\lra
H^1 (A,\OO _A)\ot \OO _{\PPP} (-q+1)\lra ...$$
$$\ \ \ \ \ \ \ \ \ \ \ \ \ \ \ \ \ \ ...
\lra H^q (A,\OO _A)\ot \OO _{\PPP} \lra 0.\ \ \ \ \ \ \ \ \
(\KK ^{\bullet}_0)$$
There is an injection of sheaves
$i:\OO _A\lra a _* \omega _X\ot P$ determined
by the choice of a section of $a_*\omega _X\ot P$.
This induces a map of complexes
$i^{\bullet}:\KK ^{\bullet}_0 \lra \KK ^{\bullet}$.

\proclaim{ Step 2}
$i_*:H^i(A,\OO _A) \lra H^i(A,a_* \omega _X \ot P)$ is an isomorphism
for $0\leq i\leq q$.
\endproclaim
We will proceed by induction. By assumption $H^0(A,\OO _A)
\cong H^0(A,a_* \omega _X\ot P)$. Assume now that $
i_*: H^i(A,\OO _A) \cong H^i(A,a_* \omega _X \ot P)$ for all $i<r$, we must
show that
$i_*:H^r(A,\OO _A) \lra H^r(A,a_* \omega _X \ot P)$ is also an isomorphism.
Twisting the complexes $\KK ^{\bullet}_0 ,\ \KK ^{\bullet}$ by $\OO _{\PPP}
(-r)$ and taking cohomology, we have:
$$
\minCDarrowwidth{.2 in}
\CD
... & & & & ...\\
@VVV & & @VVV & \\
H^{r-1}(A,\OO _A) \ot  H^{q-1}(\OO _{\PPP }(-q-1))&@>{i_*}>{\cong}>&
H^{r-1}(A,a_* \omega _X \ot P) \ot  H^{q-1}(\OO _{\PPP }(-q-1))\\
@VVV & & @VVV & \\
H^{r}(A,\OO _A)
\ot  H^{q-1}(\OO _{\PPP }(-q))&@>{i_*}>>&H^{r}
(A,a_* \omega _X \ot P) \ot  H^{q-1}(\OO _{\PPP }(-q)). \\
@VVV & & @VVV & \\
0 & & & & 0
\endCD
$$
An easy spectral sequence argument implies
that for $r>1$, the vertical lines are exact.
By the five lemma, we obtain the required isomorphism.

\proclaim{ Step 3} $a _* \omega _X \ot P =\OO _A$.
\endproclaim
Let $\QQQ$ be the cokernel of $i:\OO _A\hookrightarrow a _* \omega _X \ot
P$.
For all $i \ge 0$ and $Q\in \{ \Pic ^0(A)-\OO _A\}$, $h^i(A,\OO _A\ot Q)=h^i
(A,a _* \omega _X\ot P \ot Q)=0$. By Step 2, it follows that $h^i(A, \QQQ
\ot$ $Q ) =0$
for all $i \ge 0$ and $Q\in  \Pic ^0(A)$. By \S 1.2,
$\QQQ =0$.
This completes the proof of  (1).

(2) follows from (1) since the generic rank of $a_*\omega _X$ corresponds to
$P_1(F_{X/ \Alb (X)})$.

(3) also follows immediately from (1).

(4) By Fujita's lemma \cite {Mo, (4.1)} there is a smooth projective variety
$\tilde{X}$, and a generically
finite surjective morphism $\nu : \tilde{X} \lra X$
such that $\kappa (\tilde{X})=\kappa (X)=0$ and $P_1(\tilde{X})>0$.
Hence $P_1(\tilde{X})=1$.
The Albanese map of $\tilde{X}$ is surjective by \cite{Ka1},
and $(4)$ now follows from $(3)$.
\hfill $\square$ \enddemo

\demo{Proof of Lemma 5} Let $f_*:\Alb (X)@>>> \Alb (Y)$ be the map
induced from $f$. Since $\kappa (X)=0$, by \cite{Ka1, Theorem 1}
$\alb _X:X@>>> \Alb (X)$ is an algebraic fiber space.
It follows easily that $f_*$ and $f_* \circ \alb _X$
are surjective maps.
Consider the Stein factorization $\Alb (X)@>>> A' @>>> \Alb (Y)$
of the map $f_* $. Then $\Alb (X)@>>> A'$ is an algebraic fiber space and
$ A' @>>> \Alb (Y)$ is an \'etale map of abelian varieties.
It follows that also $a':X @>>> A'$ is an algebraic fiber space and
the fibers $F_{X/Y}$ are contracted by $a'$.
By \cite{Ka1}, there exists an induced (generically finite) rational map
$Y \dashrightarrow A'$. It follows that the generic degree of the
surjective map $Y\lra A'$ is 1.
Therefore,
$Y@>>> \Alb (Y)$ is birationally \'etale.
\hfill $\square$ \enddemo

\demo{Proof of Theorem 6}
To prove (1), we consider the generically
finite map
$\alb _Y:Y\lra \alb _Y(Y)\subset \Alb (Y)$.

\proclaim{Step 1} If $\kappa (X)\geq 0$, then
the Iitaka model of $X$ dominates the Iitaka model of $Y$.
Therefore $\kappa (X)\geq \kappa (Y)$.
\endproclaim
Let $X\lra V$ and $Y\lra W$ be appropriate birational models
of the Iitaka fibrations of $X$ and $Y$ respectively. Since $\kappa
(F_{X/V} )=0$, it follows by Lemma 5 that $f(F_{X/V})$ is birational to an
abelian
variety. And the map
 $$f(F_{X/V})\lra \alb _Y(f(F_{X/V} ))$$
is birationally
\'etale. Therefore, it is easy to see
(following the proof of \cite{Ka1, Theorem 13})
that $f(F_{X/V})$
is contained in the fibers of the Iitaka fibration $Y\lra W$.
Therefore, there exists a rational map $V\dashrightarrow W$. By changing
birational models, we may assume that it is a morphism.

\proclaim{Step 2} If $\kappa (X)=\kappa (Y)$, then $P_1(F_{X/Y})\leq 1$.
\endproclaim
Since $X\lra Y$ and $Y\lra W$ are algebraic fiber spaces and
$\dim (V)=\kappa (X)=\kappa (Y)=\dim (W)$, it
follows that $V\lra W$ is birational.
$F_{X/V}=F_{X/W}\lra F_{Y/W}$
is also an algebraic fiber space with generic fiber
$F_{X/Y}$. One sees that $\kappa (F_{X/V})=0$ and $F_{Y/W}$ is of maximal
Albanese dimension
and hence $F_{Y/W}$
is birational to an abelian variety by Lemma 5.
It follows by Theorem 4 that
$P_1(F_{X/Y})\leq 1$.

We now prove  (2).
Assume that $\kappa (Y)=0$.
Since $h^0(F_{X/Y},\omega _{X}|_{F_{X/Y}})>0$,
then $f_*\omega _{X}\ne 0$.
If $h^0(X,\omega _{X}\ot f^*P)=0$ for all $P\in \Pic ^0(Y)$, then by
Proposition 2.2,
$h^i(Y,f _*\omega _{X}\ot P)=0$ for all $i$ and all
$P\in \Pic ^0(Y)$ and hence $f_* \omega _{X}=0$ a contradiction.
Therefore $h^0(X,\omega _X \ot f^*P)>0$ for some $P\in \Pic ^0(Y)$.
By Theorem 2.1.4, we may assume that $P^{\ot r}=\OO _Y$ for some
appropriate integer $r>0$. Therefore,
$P_r(X)>0$ and $\kappa (X)\geq 0$.

If $\kappa (Y)>0$ then following \cite {Ka1, Theorem 13}, there exists
an \'etale cover $\tilde {Y}@>>> Y$ which is birational to
$\Gamma \times P$ with $\Gamma$ of general type (and maximal Albanese
dimension) and $P$ an abelian variety.
Consider the corresponding \'etale cover $\tilde {X}:=X\times _Y \tilde {Y}
@>>> X$ and the induced algebraic fiber space $\tilde {f}:\tilde {X}@>>>
\tilde {Y}$. Since $\kappa ( \tilde {X})=\kappa (X)$ and
$\kappa ( \tilde {Y})=\kappa (Y)$, it suffices to show that
$\kappa ( \tilde {X})\geq \kappa ( \tilde {Y})=\kappa (\Gamma )$.
Since $\Gamma$ is of general type,
by a theorem of Viehweg (see e.g. \cite{Mo, 6.2.d}),
we have that $\kappa (\tilde {X})\geq \kappa (\Gamma )+\kappa (F_
{\tilde {X}/\Gamma })$.
We then consider the fiber space $F_{\tilde {X}/\Gamma }
@>>> F_{\tilde {Y}/\Gamma }$. Since $\kappa
(F_{\tilde {Y}/\Gamma })=\kappa (P)=0$, one sees that $\kappa (F_{\tilde
{X}/\Gamma })
\geq \kappa (F_{\tilde {Y}/\Gamma })=0$ by the preceeding case.
 The assertion now
follows.
\hfill $\square$ \enddemo

\demo{Proof of Theorem 7} If $\dim (X)=q(X)+1$, this is \cite{Ka1, Theorem
15}.
Let $A=\Alb (X)$ and  $q:=q(X)=\dim (A)$. We have already seen that $\OO _A$
is an isolated point of $V^0(X,A,\omega _X)$.
Therefore, proceeding as in the proof of Theorem 4,
we have that
$h^q(A,a_* \omega _X) =h^0(A,a_* \omega _X)=1$.
By \cite {Ko1, Proposition 7.6}, $R^{n-q}a_* \omega _X =\omega _A$.
Therefore, by Hodge symmetry, Serre duality and \cite{Ko2}
$$h^0(X,\Omega _X^{n-q})=h^q(X,\omega _X)\geq
h^q(A,a_* \omega _X)+h^{q-(n-q)}(A,R^{n-q}a_* \omega _X) >$$
$$h^{2q-n}(A,\omega _A )= h^{n-q}(A,\OO _A)= h^0(A,\Omega ^{n-q}_A).$$
If $n=q+1$ and $P_1(X)=1$, then we would have
$h^0(X,\Omega _X^1)>h^0(A,\Omega _A^1)$
which is impossible.
\hfill $\square$
\enddemo

\proclaim {Corollary 3.4 (\cite{Ka1, Theorem 15})}
Let $\kappa (X)=0$, $\DX =q(X)+1$.
Then Conjecture K holds
\endproclaim
\demo {Proof} By Fujita's lemma \cite {Mo, (4.1)} there is
a smooth projective variety
$Y$, and a generically finite surjective morphism $\nu : Y \lra X$
such that $\kappa (Y)=\kappa (X)=0$ and $P_1(Y)=1$.
One sees that $\Alb (Y) \lra \Alb (X)$ is surjective and hence $q(Y) \ge
q(X) = \dim (X) -1$.
Therefore, by Theorem 7, $q(Y) \ge \dim (Y)$
and hence by Kawamata's Theorem,
$Y\lra \Alb (Y)$ is birational.
The fibers of $Y\lra \Alb (X) $ are translates of a fixed elliptic curve
$E$.
The corollary
now follows. \hfill $\square$ \enddemo

\enddocument
\end